\documentclass[11pt]{article}
\title{Self-Adjoint Extension of Symmetric Maps}
\author{H. N. Friedel}
\usepackage{amsmath}
\usepackage{latexsym}
\usepackage{graphicx}    % needed for including graphics e.g. EPS, PS
\begin{document}
\maketitle
\begin{abstract}
A densely-defined symmetric linear map from/to a real Hilbert space extends to a self-adjoint map.
Extension is expressed via Riesz representation.
For a case including Friedrichs extension of a strongly monotone map, self-adjoint extension is unique, and equals closure of the given map.
\end{abstract}
Let $\{ A : X \supseteq \mathcal {D} \! \mathit{o}(A) \to X\}$ be a densely-defined symmetric linear map.
Recall that if Hilbert-space $X$ is complex, then $A$ may lack self-adjoint extension (see e.g. [R]).
In contrast, self-adjoint extension must exist if our Hilbert-space is real, as will be shown here.

To prepare, we express well-known material in a form convenient for the present purpose.
For $x \in X$, let $(x | A)$ denote the linear function $\{ \mathcal {D} \! \mathit{o}(A) \ni y \mapsto (x | Ay) \}$;
we use the convention that scalar-product is linear in the second entry, conjugate-linear in the first.
Observe the adjoint domain $ \, \mathcal {D} \! \mathit{o}(A^*) \, $ equals $ \, \{ x \in X : (x | A) \  \mathrm{continuous} \} $.
Recall:  $ \mathcal {D} \! \mathit{o}(A) \, \subseteq \, \mathcal {D} \! \mathit{o}(A^*) \,$;
$ \, A$ is self-adjoint iff $\; \mathcal {D} \! \mathit{o}(A) \, = \, \mathcal {D} \! \mathit{o}(A^*) \,$.
Let $J$ denote the duality-map on $X$, which maps $x$ to function $(x | \cdot)$ in dual-space $ X^* \,$; write $J^{-1} \, = \, R$, Riesz-representation.
Extend Riesz-map $R$ so as to act on densely-defined (continuous linear) functions, such as $(x | A)$ if $\, x \in \mathcal {D} \! \mathit{o}(A^*) \,$.

\textit{Note.}  Let $A$ have symmetric extension $B$.  Then \\
(\textit{i})  $ \; \mathcal {D} \! \mathit{o}(A) \, \subseteq \, \mathcal {D} \! \mathit{o}(B) \, \subseteq \, \mathcal {D} \! \mathit{o}(B^*) \, \subseteq \, \mathcal {D} \! \mathit{o}(A^*) \; . $ \\
(\textit{ii})  $R(x | B) \, = \, R(x | A)$, if $\, x \in \mathcal {D} \! \mathit{o}(B^*) \subseteq \mathcal {D} \! \mathit{o}(A^*) \,$. \\
(\textit{iii})  $Bx \, = \, R(x|A) \,$ if $\, x \in \mathcal {D} \! \mathit{o}(B) \,$. \\
\textit{Proof.}  (\textit{i}) is known.  For $\, y \in \mathcal {D} \! \mathit{o}(A)$, see $(R(x|B) \, | \, y) = (x | By) = (x | Ay) = (R(x|A) \, | \, y) \,$;
density of $\mathcal {D} \! \mathit{o}(A)$ gives (\textit{ii}).
For $ \,x \in \mathcal {D} \! \mathit{o}(B) \,$ and $ \,y \in \mathcal {D} \! \mathit{o}(A) \,$, see $(x|A)$ is continuous, and \\
$(Bx|y) \, = \, (x|By) \, = \, (x|Ay) \, = \, ( \, R(x|A) \, | \, y) \,$; density of $\mathcal {D} \! \mathit{o}(A)$ gives (\textit{iii}).  \textit{Done.} \\
Denote by $\Lambda$ the linear map $\{ \mathcal {D} \! \mathit{o}(A^*) \ni x \mapsto R(x | A) \}$.
\textit{Note(iii)} (above) says $A$ has at-most-one symmetric extension to a given subspace $Y$,
with $\mathcal {D} \! \mathit{o}(A) \subseteq Y \subseteq \mathcal {D} \! \mathit{o}(A^*) \,$; if such extension exists, then it equals the restriction $\Lambda_{\big| Y} \;$.

\textbf{Theorem.}  Every symmetric map from/to a real Hilbert space has self-adjoint extension. \\
\textbf{Proof.}
Let $E$ denote the order-set of linear subspaces $Y$, with $\mathcal {D} \! \mathit{o}(A) \subseteq Y \subseteq \mathcal {D} \! \mathit{o}(A^*) \,$, for which restriction $\Lambda_{\big| Y} \;$ is symmetric; order by inclusion.  ($E \ni \mathcal {D} \! \mathit{o}(A)$.)
A chain C in E is bound above by the union of subspaces in $C$; so Zorn's lemma ensures $E$ has a maximal member, $Z$.
$ \, \Lambda_{\big| Z} \;$ is a maximal symmetric extension of $A$.

Write $ \, \Lambda_{\big| Z} = M\,$.
We claim $\mathcal {D} \! \mathit{o}(M) \, = \, \mathcal {D} \! \mathit{o}(M^*)$; if true, then $M$ would be self-adjoint, concluding the proof.
It is enough to show $\mathcal {D} \! \mathit{o}(M^*) \, \subseteq \, \mathcal {D} \! \mathit{o}(M)$; suppose not, seek a contradiction.
Fix $\, p \in \mathcal {D} \! \mathit{o}(M^*) \big\backslash \mathcal {D} \! \mathit{o}(M)$.
On the subspace $\, \mathcal {D} \! \mathit{o}(M) \oplus \mathbf{R}\, p \,$, define a map $T$: \\
$T(x \; + \; a \, p ) \; = \; Mx \; + \; a \, R(p|M) \ $ if $\; x \in \mathcal {D} \! \mathit{o}(M)$, $a \in \mathbf{R}$. \\
See $T$ is linear, and $T$ properly extends $M$.
To show symmetry of $T$, let $\{x,y\} \subset \mathcal {D} \! \mathit{o}(M)$ and $\{a,b\} \subset \mathbf{R}$;
note $(x \, \big| \, R(p|M)) \, = \, (p|Mx)$, $ \, (R(p|M) \, \big| \, y) \, = \, (p|My) \,$;  compute: \\
$\big( \, T(x + ap) \, \big| \, y+bp \big) \; = \; \big( \, Mx + a R(p|M) \, \big| \, y+bp \big) \; = $ \\
$(Mx|y) \; + \;  b(Mx|p) \; + \; a \big( \, R(p|M) \, \big| \, y \big) \; + \; ab \big( \, R(p|M) \, \big| \, p \big) \; = $ \\
$(x|My) \; + \; b \big( x \, \big| \, R(p|M) \, \big) \; + \; a (p|My) \; + \; ab \big( p \, \big| \, R(p|M) \, \big) \; = $ \\
$\big( x + ap \, \big| \, My + b R(p|M) \big) \; = \; \big( x + ap \, \big| \, T(y + bp) \big) \,$. \\
$M$ has symmetric proper extension $T$, so $M$ is not a maximal symmetric extension of $A$; \textit{contra}.  \textbf{Done.}

So, self-adjoint extension exists; now treat uniqueness.
Fortunately, extension is unique for some cases of interest; sometimes we may even express extension simply, as closure of the given map.
To prepare to show this, recall $A$ has symmetric closure $\bar{A} \subseteq M \,$.
Here, as before, $\{ A : X \supseteq \mathcal {D} \! \mathit{o}(A) \to X\}$ is symmetric, with self-adjoint extension $M$, from/to a Hilbert space $X$, now assumed real.
We also need the following two facts. \\
\textit{Note 1.}  If $A$ has dense image and continuous inverse, then $\bar{A}$ is the unique self-adjoint extension of $A$; $ \, M= \bar{A}$.
$\, \bar{A}$ maps onto $X$, and has continuous self-adjoint inverse. \\
\textit{Proof.}  $\bar{A}$ has dense image (since $A$ does); recall a symmetric map ($\bar{A}$) with dense image has symmetric inverse; $\bar{A}^{-1} \,$ is also closed, since $\bar{A}$ is so.
$ \, \bar{A}^{-1} \,$ equals closure of a continuous map ($ A^{-1} $), hence $\, \bar{A}^{-1} \,$ is continuous.
Since $\, \bar{A}^{-1} \,$ is closed, continuous, and has dense domain (including $\mathcal {I} \! \mathit{m}(A) \,$), we have $\mathcal {D} \! \mathit{o} \big( \bar{A}^{-1} \big)=X$.
A continuous symmetric map $\big( \bar{A}^{-1} \big)$ on the whole Hilbert space is self-adjoint.
Recall a self-adjoint map $\big( \bar{A}^{-1} \big)$ with dense image (including $\mathcal {D} \! \mathit{o}(A) \,$) has self-adjoint inverse ($\bar{A}$).
Hence $\{ \bar{A}, M \}$ are self-adjoint extensions of $A$, with $\bar{A} \subseteq M$; this forces $\bar{A}=M$, because a self-adjoint map is maximal-symmetric.  \textit{Done.} \\
\textit{Note 2.}  A (densely-defined) closed $1\!:\!1$ symmetric map has dense image. \\
\textit{Proof.}  It is enough to show $ \, p=0$, if $\, p \in \mathcal {I} \! \mathit{m} ^{\perp}(A)$ (orthogonal complement of image).
Since $\, \mathcal {D} \! \mathit{o}(A) \,$ is dense, it has a sequence $ \{ u_n \} $ converging to $p$.
If $\, x \in \mathcal {D} \! \mathit{o}(A) \, $, then \\
$ 0 \; = \; (p|Ax) \; = \; \lim (u_n |Ax) \; = \; \lim (A u_n | x). \; $
Density of $ \mathcal {D} \! \mathit{o}(A) \, $ forces $\, \lim A u_n = 0$.
$A$ is closed; ($\lim u_n = p$) and ($\lim Au_n = 0$); hence $\, p \in \mathcal {D} \! \mathit{o}(A)$, $Ap=0$.
Since $A$ is $1\!:\!1$, we have $p=0$.  \textit{Done.}

Recall (e.g. [Z]) that if our map $A$ is strongly monotone, then it has Friedrichs extension, which is self-adjoint, $1\!:\!1$, onto, with continuous self-adjoint inverse. \\
\textbf{Theorem.}  If $A$ is strongly monotone, then closure $\bar{A}$ is the unique self-adjoint extension of $A$; $\, \bar{A}$ equals Friedrichs extension. \\
\textbf{Proof.}  Let $\hat{A}$ denote Friedrichs extension; $\hat{A} \supseteq \bar{A} \,$.
Since $\hat{A}$ is $1\!:\!1$ with continuous inverse, so is its restriction $\bar{A}$.
By \textit{Note 2}, closed symmetric $1\!:\!1$ map $\bar{A}$ has dense image; then \textit{Note 1} makes $\bar{A}$ the unique self-adjoint extension of itself, and of $A$.
$\hat{A}$ is a self-adjoint extension of $A$, hence $\hat{A}=\bar{A}$.  \textbf{Done.} \\
Construction of the Friedrichs extension is complicated; how nice to express it simply (as closure), and to know it is the only self-adjoint extension.

\textbf{References} \\
$\mathrm{[R]} \ \ \ $Rudin, W.$\ \ $\textit{Functional Analysis}.  McGraw-Hill, 1991.\\
$\mathrm{[Z]} \ \ \ $Zeidler, E.$\ \ $\textit{Applied Functional Analysis: Applications to Mathematical Physics}.  Springer, 1995.
\end{document}